\documentclass[10pt]{article}
\usepackage{amssymb,amsfonts}
\usepackage{amsmath}
\usepackage{amscd,amssymb,amsthm}
\usepackage{graphics}
\usepackage{graphicx}
\usepackage{hyperref}

\graphicspath{ {./images/} }
\hypersetup{
    colorlinks=true,citecolor=blue
  }
\newtheorem{theorem}{Theorem}[section]

\newtheorem{remark}[theorem]{Remark}

\begin{document}

\title{Analysis of a map-based neuronal model}
\date{}
\author{G. Moza\thanks{
Department of Mathematics, Politehnica University of Timisoara, P-ta
Victoriei, Nr.2, 300006, Timisoara, Timis, Romania; email:
gheorghe.moza@upt.ro}, R. Efrem\thanks{
Department of Mathematics, University of Craiova, Romania} }
\maketitle

\begin{abstract}
Subthreshold oscillations in neurons are those oscillations which do not
attain the critical value of the membrane's voltage needed for triggering an
action potential (a spike). Their contribution to the forming of action
potentials in neurons is a current field of research in biology. The present
work approaches this subject using tools from mathematical modeling, more
exactly, a neuronal non-smooth map-based model is proposed and studied. The
behavior of the model in a noisy medium is also studied.
\end{abstract}

\section{Introduction}

Neuronal activity is a continuous challenge both for common people and
researchers. One of the pioneering works from this domain is \cite{hod},
where the authors describe how action potentials (electrical impulses,
spikes) are formed in the central part of a neuron. A bio-mathematical model
based on four differential equations has been proposed and studied in this
work to explain the process of forming action potentials. More other models
have been proposed since then in an effort to better understand how neurons
work. The research field is of current interest and new models based on
differential equations or difference equations (map-based) are proposed.
Neuronal map-based models have been studied lately, for example, in \cite{ba}%
--\cite{luc1}, \cite{meng}, \cite{luc4}, \cite{zm}, \cite{tig1} and in some
references therein.

While many models study various patterns of action potentials, in \cite{sil1}
is proposed a model that focuses on precursory oscillations to spiking
activity. Experimental data suggest that some neurons are brought out from
the silence regime not abruptly but slowly through small oscillations \cite%
{luc10}. These oscillations below the subthreshold of observable spikes may
shape the spiking neuronal activity when the membrane gets depolarized or
hyperpolarized \cite{luc11}. Two discrete-time phenomenological models
trying to explain how subthreshold oscillations shape spiking-bursting
neuron oscillations are proposed in \cite{sil1} and \cite{sil2}. The models
are in the form of a two-dimensional non-smooth map (discrete-time system)
given by

\begin{equation}
\begin{array}{cc}
X_{n+1}= & f_{a}\left( X_{n},Y_{n}\right) \ \ \ \ \ \ \ \ \ \ \ \ \ \  \\ 
Y_{n+1}= & Y_{n}-m\left( X_{n}+1-s\right)%
\end{array}%
,  \label{ec0a}
\end{equation}%
$n\geq 0,$ where $f_{a}:\mathbb{R}^{2}\rightarrow \mathbb{R}$ is a piecewise
differentiable function with $m\geq 0$ and $s\in 
\mathbb{R}
.$ The system is of non-smooth type because $f_{a}$ is a non-smooth
function. The $X-$variable in \eqref{ec0a} describes the dynamics of the
neuron's transmembrane potential while the parameters $a,m,s$ control
individual dynamics of the system \cite{rul1}. A recent study on two Rulkov
models of type (\ref{ec0a}) has been reported in \cite{rul17}. The map
studied in \cite{sil1} is able to generate stable subthreshold oscillations
while the map proposed in \cite{sil2} unstable subthreshold oscillations. In 
\cite{rul1} $f_{a}\left( x,y\right) $ is of the form $f_{a}\left( x,y+\beta
\right) $ to provide links with other models. In this work we aim to
introduce a new model by changing the expression of the function $f_{a}$ in (%
\ref{ec0a}). More exactly, we propose an \textit{exponential branch} in the
definition of $f_{a}$ instead of a \textit{parabola} or \textit{hyperbola}
used in the above-mentioned models \cite{sil1}, \cite{sil2}.

The paper is organized as follows. After an introductory section, in Section
2 we define our neuronal map-based two-dimensional model and point out the
first properties of it. In Sections 3 and 4 we proceed to the bifurcation
analysis of the main part of the model. In Section 5 we investigate the map
subjected to noises. Numerical studies show that the map captures the
phenomenon of vulnerability to noise of action potentials. More exactly,
subthreshold activity of the map in a noisy medium gives rise to spikes.
Moreover, we observed that our map in the presence of noise is brought out
to a spiking activity directly from the silence regime. Conclusions are
presented in the last section of the paper.

\section{Defining the model. Bifurcations}

In the present work we consider the system (\ref{ec0a}) with a function $%
f_{a}:\mathbb{R}^2\rightarrow \mathbb{R}$ given by 
\begin{equation}
f_{a}\left( x,y\right) =\left\{ 
\begin{array}{ll}
-a^{2}-e^{-a}+y, & if\ \ (x,y)\in D_1 \\ 
ax-e^{x}+y, & if\ \ (x,y)\in D_2 \\ 
a(y+1)-e^{y+1}+y, & if\ \ (x,y)\in D_3 \\ 
-1, & if\ \ (x,y)\in D_4%
\end{array}%
\right.  \label{ec0b}
\end{equation}%
where

$D_1=\{(x,y)\in\mathbb{R}^2: x<-a\},$ $D_2=\{(x,y)\in\mathbb{R}^2: x\geq
-a,\ x<y+1 \},$

$D_3=\{(x,y)\in\mathbb{R}^2: x\geq -a,\ y+1\leq x<y+2 \}$ and $D_4=\{(x,y)\in%
\mathbb{R}^2: x\geq -a,\ x\geq y+2 \}.$ 
\begin{figure}[htbp]
\begin{center}
\includegraphics[width=65mm, height=55mm]{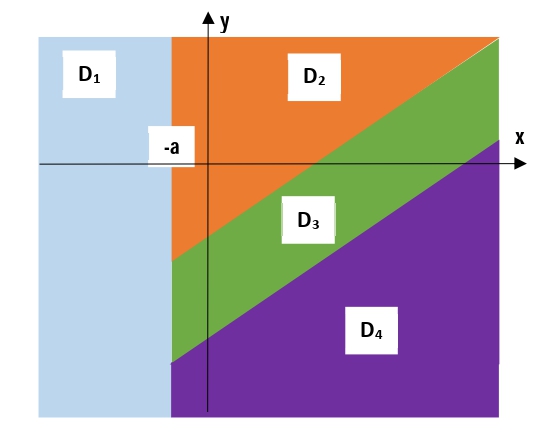}
\end{center}
\caption{The four planar sets $D_1-D_4$ forming the domain of definition of
the function $f_a$ for $a>0.$}
\label{f0}
\end{figure}

When referring further in this work at the system (\ref{ec0a}) we consider
that $f_{a}$ is given by (\ref{ec0b}). The function $f_{a}$ is piecewise
smooth. It is continuous on the borders $D_{12}$ and $D_{23}$ and
discontinuous on $D_{13},$ $D_{14}$ and $D_{34},$ where $D_{12}$ denotes the
border between the domains $D_{1}$ and $D_{2},$ and similarly for the other
borders. The domains $D_1-D_4$ and their borders are sketched in Fig.\ref{f0}. For those values of $m$ small enough, $0<m\ll 1,$ the system (\ref%
{ec0a}) is of slow-fast type; the first equation in \eqref{ec0a} is termed
the fast subsystem while the second the slow subsystem.

\begin{remark}
The reason we propose an exponential branch in the expression of the
function $f_{a}$ is to see if exponential terms bring new types of
oscillations of the map (\ref{ec0a}), new patterns of oscillations both in
the subthreshold and spiking activity regimes, or a larger spectrum of
stability, compared with the models studied in \cite{sil1} and \cite{sil2}
which use \textit{parabola} or \textit{hyperbola} terms.
\end{remark}

The main dynamics of the system (\ref{ec0a}) is described by the second
branch of $f_{a}$ since it has the both variables $x$ and $y.$ The other
branches are some kind of limiters. Hence, we will study in the following
the system

\begin{equation}
\begin{array}{cc}
X_{n+1}= & aX_{n}-e^{X_{n}}+Y_{n} \\ 
Y_{n+1}= & Y_{n}-m\left( X_{n}+1-s\right)%
\end{array}
\label{ec0}
\end{equation}%
with $-a\leq X_{n}<Y_{n}+1,$ and $a,s\in \mathbb{R},$ $m\geq 0.$

We study first (\ref{ec0}) at $m=0.$ In this case, the second equation in (%
\ref{ec0}) is independent from the first and $Y_{n}=Y_{0},$ for all $n\geq
0. $ Hence, $Y_{0}$ becomes a real parameter in the first equation.
Depending on $a$ and $Y_{0},$ the system can have zero, one or two fixed
points. More exactly, if $a<1$ it has a single fixed point for any fixed
real $Y_{0},$ at $a=1$ it has no fixed points for $Y_{0}\leq 0$ and one
fixed point $X=\ln Y_{0}$ for $Y_{0}>0,$ while for $a>1$ it can have zero,
one or two fixed points. Let us detail the latter case, $a-1>0.$ Denote by 
\begin{equation*}
h\left( x\right) =\left( a-1\right) x-e^{x}+Y_{0}
\end{equation*}%
and $x_{0}=\ln \left( a-1\right) $ the root of $h^{\prime }\left( x\right)
=0.$ Notice that $h\left( \infty \right) =-\infty $ and $h\left( -\infty
\right) =-\infty .$ Therefore, if $h\left( x_{0}\right) <0$ the system (\ref%
{ec0}) has no fixed points while if $h\left( x_{0}\right) =0,$ it has a
unique fixed point $X=x_{0}.$ When $h\left( x_{0}\right) >0,$ the system has
two different fixed points denoted by $X_{1,2}\left( a,Y_{0}\right) $ which
satisfy $X_{1}<\ln \left( a-1\right) <X_{2}.$ From 
\begin{equation*}
\frac{\partial X_{n+1}}{\partial X_{n}}=a-e^{X_{n}},
\end{equation*}%
it follows that $X_{1}$ is unstable whenever it exists, that is on $X<\ln
\left( a-1\right) ,$ respectively, $X_{2}$ is stable for $\ln
(a-1)<X_{2}<\ln (a+1)$ and unstable for $X_{2}>\ln (a+1).$ For example, if $%
a=e+1$ and $Y_{0}=1,$ we find $X_{1}=0$ and $X_{2}\simeq 1.752.$

Assume further $m>0.$ The system (\ref{ec0}) has a unique fixed point $%
A\left( s-1,\left( 1-a\right) \left( s-1\right) +e^{s-1}\right) $ in this
case, provided that the parameters satisfy the inequality $%
-a<s-1<\allowbreak a+\left( 1-a\right) s+e^{s-1}.$ Introducing the new
coordinates

\begin{equation*}
\left( 
\begin{array}{c}
X \\ 
Y%
\end{array}%
\right) =\left( 
\begin{array}{c}
x+s-1 \\ 
y+\left( 1-a\right) \left( s-1\right) +e^{s-1}%
\end{array}%
\right) ,
\end{equation*}

\noindent the map (\ref{ec0}) becomes

\begin{equation}
\begin{array}{cc}
\tilde{x}= & y+ax+e^{s-1}-e^{x+s-1} \\ 
\tilde{y}= & y-mx%
\end{array}
\label{ec1}
\end{equation}

\noindent where $a,s\in \mathbb{R}$ and $m>0.$ In (\ref{ec1}) and in other
discrete systems below we will use the notations $\left( \tilde{x},\tilde{y}%
\right) =\left( x_{n+1},y_{n+1}\right) $ and $\left( x,y\right) =\left(
x_{n},y_{n}\right) $ for all $n\geq 0.$ The new map (\ref{ec1}) is defined
on 
\begin{equation*}
-a\leq s-1+x<a+\left( 1-a\right) s+e^{s-1}+y
\end{equation*}%
and has a unique fixed point $O\left( 0,0\right) .$ Its Jacobian at $O$ is $%
J=\left( 
\begin{array}{cc}
a-e^{s-1} & 1 \\ 
-m & 1%
\end{array}%
\right) $ with the eigenvalues (multipliers)

\begin{equation}
\lambda ^{\pm }=\frac{1}{2}\left( a-e^{s-1}+1\pm \sqrt{\left(
e^{s-1}-a+1\right) ^{2}-4m}\right) .  \label{ec1a}
\end{equation}%
Since $m>0,$ two main bifurcations may appear in (\ref{ec1}), namely
period-doubling and Neimark-Sacker.

\section{Period-doubling bifurcations}

Assume first $\left( e^{s-1}-a+1\right) ^{2}\geq 4m,$ hence $\lambda ^{\pm }$
are real. Notice that $\lambda ^{\pm }\neq 1$ since $m>0,$ thus, a fold
bifurcation is not possible. Denote further by 
\begin{equation*}
m_{0}=2\left( e^{s-1}-a-1\right)
\end{equation*}
whenever $e^{s-1}-a-1>0.$ Then, $\lambda ^{+}=-1$ and $\lambda
^{-}=\allowbreak a-e^{s-1}+2$ if $m_{0}>4$ and $m=m_{0},$ respectively, $%
\lambda ^{-}=-1$ and $\lambda ^{+}=\allowbreak a-e^{s-1}+2$ if $0<m_{0}<4$
and $m=m_{0}.$ At $m=m_{0}=4,$ we have $\lambda ^{+}=\lambda ^{-}=-1.$ Thus,
a \textit{period-doubling bifurcation} occurs in the system (\ref{ec1}) at $%
m=m_{0},$ which is also known as the \textit{flip} bifurcation. In order to
study this bifurcation we need to reduce the system (\ref{ec1}) to the
center manifold through $O\left( 0,0\right) .$ The next theorem from \cite%
{kuz} will be useful to our study.

\begin{theorem}
\label{t2} Assume the one-dimensional discrete system 
\begin{equation}
\tilde{x}=f(x,m),  \label{ec4}
\end{equation}%
$f$ sufficiently smooth, $x,\tilde{x},m\in 
\mathbb{R}
,$ has at $m=0$ a fixed non-hyperbolic point $x_{0}=0$ with the multiplier $%
\mu =f_{x}(0,0)=-1.$ Consider that the following two nondegeneracy (generic)
conditions hold:\newline
pd1) $c\left( 0\right) \overset{def}{=}\frac{1}{4}\left( \frac{\partial ^{2}f%
}{\partial x^{2}}(0,0)\right) ^{2}+\frac{1}{6}\frac{\partial ^{3}f}{\partial
x^{3}}(0,0)\neq 0,$\newline
pd2) $\frac{\partial ^{2}f}{\partial x\partial m}(0,0)\neq 0.$\newline
Then, the system \eqref{ec4} is locally topologically equivalent near the
origin to the following normal form 
\begin{equation}
\tilde{\xi}=-(1+\beta \left( m\right) )\xi +s\xi ^{3},  \label{ec4a}
\end{equation}%
where $s=sign\left( c\left( 0\right) \right) $\ and $\beta \left( m\right)
=-1-\frac{\partial f}{\partial x}\left( 0,m\right) .$
\end{theorem}

To this end, consider further the extended form of (\ref{ec1}) 
\begin{eqnarray}
\tilde{m} &=&m  \notag \\
\tilde{x} &=&y+x\left( a-e^{s-1}\right) -\frac{1}{2}x^{2}e^{s-1}-\frac{1}{6}%
x^{3}e^{s-1}+\allowbreak O\left( x^{4}\right)  \label{ec3} \\
\tilde{y} &=&y-mx  \notag
\end{eqnarray}%
which has the non-hyperbolic fixed point $\left( m_{0},0,0\right) $ with the
multipliers $\lambda _{1}=1,$ $\lambda _{2}=-1$ and $\lambda
_{3}=a-e^{s-1}+2.$ Assume further $m_{0}\neq 4,$ thus, $\lambda _{3}\neq \pm
1.$ We used the Taylor expansion of $e^{x+s-1}$ at $x=0$ in (\ref{ec3}).
Hence, there exists a 2-dimensional center manifold $W^{c}$ of equation $%
y=y\left( m,x\right) $ and $\left( m_{0},0,0\right) \in W^{c},$ given
locally by: 
\begin{align}
W^{c}:y& =a_{1}\left( m-m_{0}\right) +a_{2}\left( m-m_{0}\right)
^{2}+a_{3}\left( m-m_{0}\right) ^{3}+b_{1}x+b_{2}x^{2}+b_{3}x^{3}
\label{ec3a} \\
& +c_{1}\left( m-m_{0}\right) x+c_{2}\left( m-m_{0}\right) x^{2}+c_{3}\left(
m-m_{0}\right) ^{2}x+...  \notag
\end{align}%
Then 
\begin{align*}
\tilde{y}& =a_{1}\left( \tilde{m}-m_{0}\right) +a_{2}\left( \tilde{m}%
-m_{0}\right) ^{2}+a_{3}\left( \tilde{m}-m_{0}\right) ^{3}+b_{1}\tilde{x}%
+b_{2}\tilde{x}^{2}+b_{3}\tilde{x}^{3} \\
& +c_{1}\left( \tilde{m}-m_{0}\right) \tilde{x}+c_{2}\left( \tilde{m}%
-m_{0}\right) \tilde{x}^{2}+c_{3}\left( \tilde{m}-m_{0}\right) ^{2}\tilde{x}%
+...
\end{align*}

Using now \eqref{ec3} and \eqref{ec3a}, we find that $%
a_{1}=0,a_{2}=0,a_{3}=0,$ $b_{1}=\frac{m_{0}}{2},\ \ b_{2}=\frac{1}{2}%
e^{s-1},\ \ b_{3}=-\frac{1}{6}e^{s-1}\frac{m_{0}}{4-m_{0}},$ $c_{1}=\frac{2}{%
4-m_{0}},\ \ c_{2}=\frac{4e^{s-1}}{m_{0}\left( 4-m_{0}\right) },\ \ c_{3}=%
\frac{8}{\left( 4-m_{0}\right) ^{3}}.$

Finally, the system (\ref{ec3}) restricted to the central manifold is
1-dimensional and for $\left\vert m-m_{0}\right\vert $ small enough it reads

\begin{equation}
\tilde{x}=x\sigma _{1}+x^{2}\sigma _{2}+x^{3}\sigma _{3}+O\left(
x^{4}\right) \overset{not}{=}f\left( x,m\right)  \label{ec4b}
\end{equation}

\noindent where $\sigma _{1}=\frac{2\left( 4-m_{0}\right) ^{2}\left(
m-m_{0}\right) +8\left( m-m_{0}\right) ^{2}}{\left( 4-m_{0}\right) ^{3}}-1,$ 
$\sigma _{2}=\frac{4e^{s-1}\left( m-m_{0}\right) }{m_{0}\left(
4-m_{0}\right) }$ and $\sigma _{3}=-\frac{2e^{s-1}}{3\left( 4-m_{0}\right) }%
. $

Because the hyperplanes $U_{m}=\left\{ \left( m,x,y\right) :m=const\right\} $
are invariant with respect to the map (\ref{ec3}), the manifold $W^{c}$ is
foliated by the 1-dimensional invariant manifolds 
\begin{equation*}
W_{m}^{c}=W^{c}\cap U_{m}.
\end{equation*}%
This manifold $W_{m}^{c}$ is a local invariant parameter-dependent manifold
of the initial system (\ref{ec1}). Then (\ref{ec4b}) has the multiplier $\mu
=\frac{\partial f}{\partial x}(0,m_{0})=-1$ while $\frac{\partial ^{2}f}{%
\partial x\partial m}(0,m_{0})=\frac{\partial \sigma _{1}}{\partial m}\left(
m_{0}\right) =-\frac{2}{m_{0}-4}\neq 0;$ $m_{0}\neq 4.$ Further, 
\begin{equation*}
c\left( 0\right) =\frac{1}{4}\left( \frac{\partial ^{2}f}{\partial x^{2}}%
\right) ^{2}+\frac{1}{6}\frac{\partial ^{3}f}{\partial x^{3}}=\sigma
_{2}^{2}+\sigma _{3}=\frac{2e^{s-1}}{3\left( m_{0}-4\right) }\neq 0.
\end{equation*}%
Hence pd1)-2) are satisfied. The coefficient $c\left( 0\right) $ can be
obtained also by another method \cite{kuz}. Indeed, denote by $z=\left( 
\begin{array}{c}
x \\ 
y%
\end{array}%
\right) .$ Using the Taylor expansion of $e^{x+s-1}$ at $x=0,$ the system (%
\ref{ec1}) when $m=m_{0}$ can be put in the form

\begin{equation*}
\tilde{z}=Jz+\frac{1}{2}B\left( z,z\right) +\frac{1}{6}C\left( z,z,z\right) ,
\end{equation*}%
where $B\left( z,u\right) =\left( 
\begin{array}{c}
-xx_{1}e^{s-1} \\ 
0%
\end{array}%
\right) ,$ $C\left( z,u,v\right) =\left( 
\begin{array}{c}
-xx_{1}x_{2}e^{s-1} \\ 
0%
\end{array}%
\right) $ and $u=\left( 
\begin{array}{c}
x_{1} \\ 
y_{1}%
\end{array}%
\right) ,$ $v=\left( 
\begin{array}{c}
x_{2} \\ 
y_{2}%
\end{array}%
\right) $ are two vectors. An eigenvector for $\lambda ^{-}=-1$ is $q=\left( 
\begin{array}{c}
1 \\ 
e^{s-1}-a-1%
\end{array}%
\right) .$ An adjoint eigenvector $p$ to $q,$ i.e. $J^{T}p=-p,$ such that $%
\left\langle p,q\right\rangle =1,$ is $p=\allowbreak \left( 
\begin{array}{c}
\frac{2}{a-e^{s-1}+3} \\ 
-\frac{1}{a-e^{s-1}+3}%
\end{array}%
\right) .$ We can determine now the coefficient $c\left( 0\right) $ from the
formula 
\begin{eqnarray*}
c\left( 0\right) &=&\frac{1}{6}\left\langle p,C\left( q,q,q\right)
\right\rangle -\frac{1}{2}\left\langle p,B\left( q,\left( J-I_{2}\right)
^{-1}B\left( q,q\right) \right) \right\rangle \\
&=&\frac{1}{6}\left\langle p,C\left( q,q,q\right) \right\rangle =\frac{%
4e^{s-1}}{m_{0}-4}.
\end{eqnarray*}

We notice the expressions of $c\left( 0\right) $ obtained with the two
methods differ by a positive constant but $s=sign\left( c\left( 0\right)
\right) $ is the same. In fact $c\left( 0\right) $ is not unique, it depends
on the vectors $p$ and $q.$

Hence, a period-doubling bifurcation occurs at $m=m_{0}$ and the system %
\eqref{ec4b} is locally topologically equivalent near the origin to the
normal form 
\begin{equation}
\tilde{\xi}=-(1+\beta \left( m\right) )\xi +s\xi ^{3},  \label{ec5}
\end{equation}%
where $s=+1$ if $m_{0}>4$ and $s=-1$ if $0<m_{0}<4.$ Here 
\begin{equation*}
\beta \left( m\right) =-1-\sigma _{1}=\allowbreak 2\left( 4\left(
m-m_{0}\right) +\left( m_{0}-4\right) ^{2}\right) \frac{m-m_{0}}{\left(
m_{0}-4\right) ^{3}}
\end{equation*}%
with $\beta \left( m_{0}\right) =0$ and $\frac{\partial \beta }{\partial m}%
\left( m_{0}\right) =\allowbreak \frac{2}{m_{0}-4}\neq 0.$

This means that a cycle of period two bifurcates at $m=m_{0}$ in (\ref{ec5}%
). Assume first $s=-1$ and $0<m_{0}<4.$ Denoting by $g\left( \xi \right)
=-(1+\beta )\xi -\xi ^{3},$ one can show that $g^{2}\left( \xi \right)
=g\left( g\left( \xi \right) \right) $ has two nontrivial fixed points
depending on a sufficiently small $\beta $ with $\beta <0$ given by 
\begin{equation*}
\xi _{1,2}=\pm \sqrt{-\beta }+O\left( \beta \right) ,
\end{equation*}%
that is, $g_{m}^{2}\left( \xi _{1,2}\right) =\xi _{1,2}.$ These two points
form a period-two \textit{unstable} cycle of the map (\ref{ec4b}), i.e. $\xi
_{1}=f\left( \xi _{2},m\right) $ and $\xi _{2}=f\left( \xi _{1},m\right) .$
The cycle vanishes for $\beta \geq 0.$ The higher-order terms in (\ref{ec4b}%
) do not affect the bifurcation. In the second case, $s=+1$ and $m_{0}>4,$ $%
g^{2}\left( \xi \right) $ has two nontrivial fixed points in the form $\xi
_{1,2}=\pm \sqrt{\beta }+O\left( \beta \right) $ which form a \textit{stable 
}cycle of the map (\ref{ec4b}). Concluding, we have the following theorem.

\begin{theorem}
\label{t1} If $0<m\leq \frac{\left( m_{0}+4\right) ^{2}}{16}$ the system (%
\ref{ec1}) undergoes a period-doubling bifurcation at $m_{0}.$ More exactly,
if $0<m_{0}<4$ (respectively, $m_{0}>4$ ) a period-two unstable
(respectively, stable) cycle is born for $0<m<m_{0}-\frac{1}{4}\left(
m_{0}-4\right) ^{2}$ or $m>m_{0}.$
\end{theorem}

\begin{remark}
Of practical interest is the case $0\leq m\ll 1$ and, because the
bifurcation occurs at $m=m_{0},$ this yields $0\leq m_{0}\ll 1.$ Thus, we
are in the case with $s=-1$ and $\beta \left( m\right) <0,$ which is
equivalent to $m>m_{0},$ since $m_{0}$ is sufficiently small. It follows
that, the period-doubling bifurcations may occur when the parameters lie in
the zones of practical interest and contribute to the formation of action
potentials in the model (\ref{ec0a}). This represents an added value of our
model compared to other models (see e.g. \cite{sil1}) in which
period-doubling bifurcations occur outside the zone of practical interest.
\end{remark}

\section{Neimark-Sacker bifurcation}

In order to study this bifurcation in our map (\ref{ec1}), we describe it in
a general setup first to ease its presentation (see e.g. \cite{kuz}).
Consider a discrete-time system 
\begin{equation*}
\tilde{x}=f(x,\alpha ),x=(x_{1},x_{2})\in \mathbb{R}^{2},\alpha \in \mathbb{R%
},
\end{equation*}%
$f$ sufficiently smooth, having for all $|\alpha |$ small enough a fixed
point $x_{0}=\left( 0,0\right) $ with the multipliers $\mu _{1,2}(\alpha
)=r(\alpha )e^{\pm i\theta (\alpha )},$ where $r(0)=1,$ $\theta (0)=\theta
_{0}$ and $0<\theta _{0}<\pi .$ The system can be put in the form 
\begin{equation}
\tilde{x}=A(\alpha )x+F(x,\alpha ),  \label{nm1}
\end{equation}%
where $\mu _{1,2}(\alpha )$ are the eigenvalues of the Jacobian matrix $%
A(\alpha ).$ The radius $r$ can also be written in the form $r(\alpha
)=1+\beta (\alpha )$ for some smooth functions $\beta \left( \alpha \right) $
with $\beta (0)=0.$ Assume $\beta ^{\prime }(0)\neq 0$ and consider $\beta $
as a new parameter of the system. Then the multipliers are of the form $\mu
(\beta ),$ $\overline{\mu }(\beta )$ where $\mu (\beta )=(1+\beta
)e^{i\theta (\beta )},$ with $\theta (\beta )$ a smooth function such that $%
\theta (0)=\theta _{0}.$ In complex variables $z=x_{1}+ix_{2},$ $\bar{z}%
=x_{1}-ix_{2},$ the equation (\ref{nm1}) can be put, for $|\beta |$ small
enough, in the form 
\begin{equation*}
\tilde{z}=\mu (\beta )z+g(z,\overline{z},\beta ),
\end{equation*}%
where $\beta \in \mathbb{R},z\in \mathbb{C},$ $\mu (\beta )=(1+\beta
)e^{i\theta (\beta )},$ $\mu _{0}=e^{i\theta _{0}},$ and $g$ is a smooth
complex function of $z,\overline{z},\beta $ whose Taylor expansion in $(z,%
\overline{z})$ begins with at least quadratic terms, 
\begin{equation*}
g(z,\overline{z},\beta )=\sum_{i+j\geq 2}\frac{1}{i!j!}g_{ij}(\beta )z^{i}%
\overline{z}^{j}.
\end{equation*}

\begin{theorem}
\label{t5}(Generic Neimark-Sacker bifurcation) Assume the two-dimensional
one-parameter system 
\begin{equation}
\tilde{x}=f(x,\alpha )  \label{tns}
\end{equation}%
$x\in 
\mathbb{R}
^{2},$ $\alpha \in 
\mathbb{R}
,$ has for all $|\alpha |$ small enough a fixed point $x_{0}=\left(
0,0\right) $ with the multipliers $\mu _{1,2}(\alpha )=r(\alpha )e^{\pm
i\theta (\alpha )},$ where $r(0)=1,$ $\theta (0)=\theta _{0},$ such that $%
\frac{dr}{d\alpha }(0)\neq 0$ and $e^{ik\theta _{0}}\neq 1$ for $k=1,2,3,4.$
Denote by 
\begin{equation}
c_{1}(0)=\frac{g_{20}\left( 0\right) g_{11}\left( 0\right) \left( 1-2\mu
_{0}\right) }{2\left( \mu _{0}^{2}-\mu _{0}\right) }+\frac{\left\vert
g_{11}\left( 0\right) \right\vert ^{2}}{1-\bar{\mu}_{0}}+\frac{\left\vert
g_{02}\left( 0\right) \right\vert ^{2}}{2\left( \mu _{0}^{2}-\bar{\mu}%
_{0}\right) }+\frac{g_{21}\left( 0\right) }{2},  \label{ec12c}
\end{equation}%
where $\mu _{0}=e^{i\theta _{0}}.$ Then, the map (\ref{tns}) is locally
topologically equivalent near the origin $x_{0}$ for all $|\alpha |$ small
enough to the map in polar coordinates given by 
\begin{equation}
\left\{ 
\begin{array}{cc}
\tilde{\rho}= & \rho \left( 1+\beta +d\left( \beta \right) \rho ^{2}\right)
+\rho ^{4}R\left( \beta ,\rho \right) \\ 
\tilde{\varphi}= & \varphi +\theta \left( \beta \right) +\rho ^{2}Q\left(
\beta ,\rho \right)%
\end{array}%
\right.  \label{pol}
\end{equation}%
where $\beta \left( \alpha \right) =r\left( \alpha \right) -1$ and $d\left(
0\right) =Re\left( e^{-i\theta _{0}}c_{1}(0)\right) .$
\end{theorem}

The first equation in (\ref{tns}), the $\rho -$map, is independent of $%
\varphi $ and can be analyzed separately as a one-dimensional map. It has
two fixed points, $\rho _{0}=0$ for all $\beta ,$ respectively, $\rho _{1}=%
\sqrt{\frac{-\beta }{d\left( \beta \right) }}+O\left( \beta \right) ,$
defined for $d\left( 0\right) \neq 0$ and $\beta d\left( \beta \right) <0,$ with $|\beta|$ sufficiently small.
It follows from the $\rho -$map that $\rho _{0}=0$ is (linearly) stable if $%
\beta <0$ and (linearly) unstable if $\beta >0.$ At $\beta =0,$ the
stability of $\rho _{0}=0$ depends on $d\left( 0\right) ,$ more exactly, it
is (nonlinearly) stable if $d\left( 0\right) <0,$ respectively,
(nonlinearly) unstable if $d\left( 0\right) >0.$

When $d\left( 0\right) <0,$ the second fixed point $\rho _{1}$ exists and is
stable for $\beta >0,$ while if $d\left( 0\right) >0,$ $\rho _{1}$ exists
and is unstable for $\beta <0.$ To the fixed point $\rho _{1}$ in the $\rho
- $map, it corresponds a closed invariant curve in the system (\ref{tns}).

\begin{remark}
\label{rem1a}If $d(0)\neq 0,$ there exists a neighborhood of $x_{0}$ in
which a unique closed invariant curve $\Gamma _{\beta }$ bifurcates from $%
x_{0}$ as $\beta $ crosses $0.$ More exactly, if $d(0)<0,$ the curve $\Gamma
_{\beta }$ exists for $\beta >0$ and is unique and stable. It vanishes for $%
\beta \leq 0.$ Similarly, if $d(0)>0,$ the curve $\Gamma _{\beta }$ exists
for any $\beta <0,$ being unique and unstable. It disappears for $\beta \geq
0.$
\end{remark}

The Neimark-Sacker (NS) bifurcation is known also as the Andronov-Hopf
bifurcation (or simply Hopf for discrete systems). We study now the NS
bifurcation in the map (\ref{ec1}). To this end, assume that the above
multipliers $\lambda ^{\pm }$ are strictly complex, that is $4m-\left(
e^{s-1}-a+1\right) ^{2}>0$ or $m>\frac{\left( m_{0}+4\right) ^{2}}{16}\geq
0. $ Consider $a\in 
\mathbb{R}
$ as the bifurcation parameter. Then $\lambda ^{\pm }\left( a\right)
=h\left( a\right) \pm i\omega \left( a\right) $ with $h\left( a\right) =%
\frac{1}{2}\left( a-e^{s-1}+1\right) $ and $\omega \left( a\right) =\frac{1}{%
2}\sqrt{4m-\left( 1+e^{s-1}-a\right) ^{2}}.$ In exponential form, they read 
\begin{equation*}
\lambda ^{\pm }=r\left( a\right) e^{\pm i\theta \left( a\right) }
\end{equation*}%
where $r\left( a\right) =\allowbreak \allowbreak \sqrt{a+m-e^{s-1}}$ and $%
\tan \left( \theta \left( a\right) \right) =\frac{\omega \left( a\right) }{%
h\left( a\right) },$ with $a+m-e^{s-1}\geq 0.$ Two complex eigenvectors
corresponding to the eigenvalues $\lambda ^{\pm }(a)$ are

\begin{equation*}
u^{\pm }=\left( 
\begin{array}{cc}
\omega \left( a\right) \pm i\left( 1-h\left( a\right) \right) & \pm mi%
\end{array}
\right)^{T}
\end{equation*}%
($T$ stands for transpose). It implies that $p=\left( u^{+}+u^{-}\right)
/2=\left( 
\begin{array}{cc}
\omega \left( a\right) & 0%
\end{array}%
\right) ^{T}$ and $q=\left( u^{+}-u^{-}\right) /2i=\left( 
\begin{array}{cc}
1-h\left( a\right) & m%
\end{array}%
\right) ^{T}$ are two real independent vectors which lead to the
transformation

\begin{equation}
\left( 
\begin{array}{c}
x \\ 
y%
\end{array}%
\right) =\left( 
\begin{array}{cc}
1-h\left( a\right) & \omega \left( a\right) \\ 
m & 0%
\end{array}%
\right) \left( 
\begin{array}{c}
u \\ 
v%
\end{array}%
\right) .  \label{ec13a}
\end{equation}

Using Taylor expansion of the function $e^{x+s-1}$ at $x=0,$ the system %
\eqref{ec1} reads

\begin{equation*}
\left( 
\begin{array}{c}
\tilde{x} \\ 
\tilde{y}%
\end{array}%
\right) =\left( 
\begin{array}{cc}
a-e^{s-1} & 1 \\ 
-m & 1%
\end{array}%
\right) \left( 
\begin{array}{c}
x \\ 
y%
\end{array}%
\right) +\left( 
\begin{array}{c}
-\frac{x^{2}}{2}e^{s-1}-\frac{x^{3}}{6}e^{s-1} \\ 
0%
\end{array}%
\right) +...
\end{equation*}

\noindent which, by using \eqref{ec13a}, leads to 
\begin{equation}
\left( 
\begin{array}{c}
\tilde{u} \\ 
\tilde{v}%
\end{array}%
\right) =\left( 
\begin{array}{cc}
h\left( a\right) & -\omega \left( a\right) \\ 
\omega \left( a\right) & h\left( a\right)%
\end{array}%
\right) \left( 
\begin{array}{c}
u \\ 
v%
\end{array}%
\right) -\left( 
\begin{array}{c}
0 \\ 
\frac{1}{6\omega }e^{s-1}\left( u+v\omega -hu\right) ^{2}\left( u+v\omega
-hu+3\right)%
\end{array}%
\right) \allowbreak +...  \label{ec14}
\end{equation}

\noindent Finally, we can write the system \eqref{ec14} in complex variables 
$z=u+iv,$ $\bar{z}=u-iv,$ in the form

\begin{equation*}
\tilde{z}=r\left( a\right) e^{i\theta \left( a\right) }z+\frac{1}{2}%
g_{20}z^{2}+g_{11}z\bar{z}+\frac{1}{2}g_{02}\bar{z}^{2}+\frac{1}{2}%
g_{21}z^{2}\bar{z}+...
\end{equation*}

\noindent where the coefficients are $g_{20}=\frac{e^{s-1}}{4}\frac{\left(
-ih+\omega +i\right) ^{2}i}{\omega },\ g_{11}=-\frac{ie^{s-1}}{4}\frac{%
\left( h-1\right) ^{2}+\omega ^{2}}{\omega },$ $g_{02}=\frac{ie^{s-1}}{4}%
\frac{\left( ih+\omega -i\right) ^{2}}{\omega }$ and $g_{21}=\frac{e^{s-1}}{%
8\omega }\left( ih+\omega -i\right) \left( h+i\omega -1\right) ^{2},$ with $%
h=h\left( a\right) $ and $\omega =\omega \left( a\right) .$

The first condition leading to the NS bifurcation is $r\left( a\right) =1,$
which yields the surface 
\begin{equation*}
(NS):a=e^{s-1}-m+1\overset{not}{=}a_{NS}.
\end{equation*}

On this surface, $\tan \theta _{0}=\frac{\sqrt{4m-m^{2}}}{2-m},$ $m\in
\left( 0,4\right) ,$ $m\neq 2,$ $h_{0}=\allowbreak \frac{1}{2}\left(
2-m\right) ,$ $\omega _{0}=\allowbreak \frac{1}{2}\sqrt{4m-m^{2}},$ $%
e^{i\theta _{0}}=\lambda ^{+}\left( a_{NS}\right) =h_{0}+i\omega _{0}$ with $%
0<\theta _{0}<\pi .$ The first nondegeneracy condition $\frac{\partial r}{%
\partial a}\left( a_{NS}\right) =\frac{1}{2}\neq 0$ holds true independent
on the parameters, while the second $e^{ik\theta _{0}}\neq 1$ is satisfied
provided that $m\in \left( 0,4\right) \smallsetminus \left\{ 2,3\right\} .$

Expressing the coefficients $g_{ij}$ at\ $a=a_{NS}$ and using \eqref{ec12c}
we obtain

\begin{equation}
d(0)=-\frac{e^{s-1}}{16}m\left( 1+e^{s-1}\right) \allowbreak .  \label{ec14a}
\end{equation}%
A similar result for $a(0)$ can be obtained by another method (see \cite{kuz}
, pag. 186), namely 
\begin{eqnarray*}
d(0) &=&\frac{1}{2}Re\left\{ e^{-i\theta _{0}}\left[ 2\left\langle p,B\left(
q,\left( I-J\right) ^{-1}B\left( q,\bar{q}\right) \right) \right\rangle
+\left\langle p,B\left( \bar{q},\left( e^{2i\theta _{0}}I-J\right)
^{-1}B\left( q,q\right) \right) \right\rangle \right] \right\} \\
&&+\frac{1}{2}Re\left\{ e^{-i\theta _{0}}\left\langle p,C\left( q,q,\bar{q}%
\right) \right\rangle \right\} \\
&=&-\frac{e^{s-1}}{16m\left( 4-m\right) }\left( 1+e^{s-1}\right) ,
\end{eqnarray*}

\noindent where the vectors $p=\allowbreak \left( 
\begin{array}{c}
m+i\sqrt{4m-m^{2}} \\ 
-2%
\end{array}%
\right) $ and $q=\allowbreak \allowbreak \left( 
\begin{array}{c}
\frac{1}{2}\frac{i}{\sqrt{4m-m^{2}}} \\ 
-\frac{1}{4}\frac{m+i\sqrt{4m-m^{2}}-4}{m-4}%
\end{array}%
\right) $ satisfy $Jq=\lambda ^{+}q,$ $J^{T}p=\lambda ^{-}p$ and $%
\left\langle p,q\right\rangle =1;$ $B$ and $C$ are given in the previous
section. Since $m\left( 4-m\right) >0,$ the two expressions of $d(0)$
coincide up to a positive constant. Only their sign matters in the
following. Concluding, we have the following theorem.

\begin{theorem}
\label{t3} If $\frac{\left( m_{0}+4\right) ^{2}}{16}<m<4$ and $m\notin
\left\{ 2,3\right\} $ the system (\ref{ec1}) undergoes a NS bifurcation at $%
a_{NS}.$ There exists a neighborhood of $x_{0}=\left( 0,0\right) $ in which
a unique closed invariant curve $\Gamma _{a}$ bifurcates from $x_{0}$ as $a$
crosses $a_{NS}.$ More exactly, $\Gamma _{a}$ is stable (by Remark \ref%
{rem1a}) and exists for all $s\in 
\mathbb{R}
$ and $a>a_{NS},$ with $\left\vert a-a_{NS}\right\vert $ sufficiently small.
\end{theorem}

\begin{figure}[h]
\begin{center}
\includegraphics[width=55mm, height=55mm]{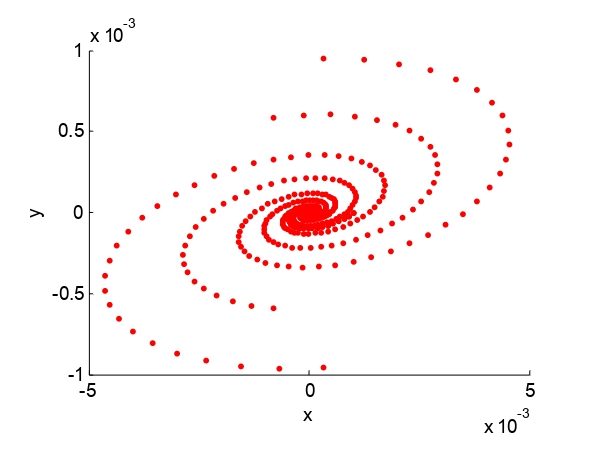} %
\includegraphics[width=55mm, height=55mm]{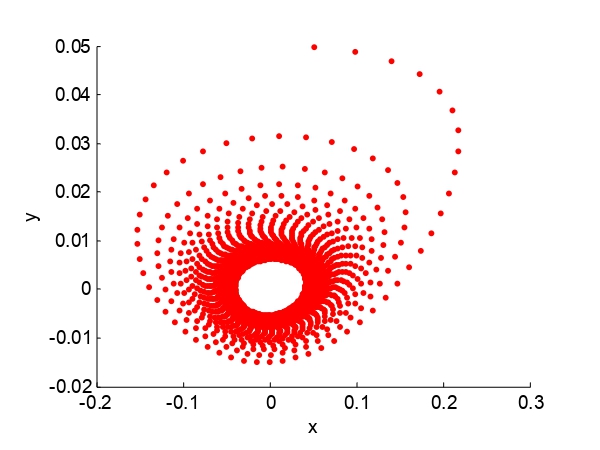} %
\includegraphics[width=55mm, height=55mm]{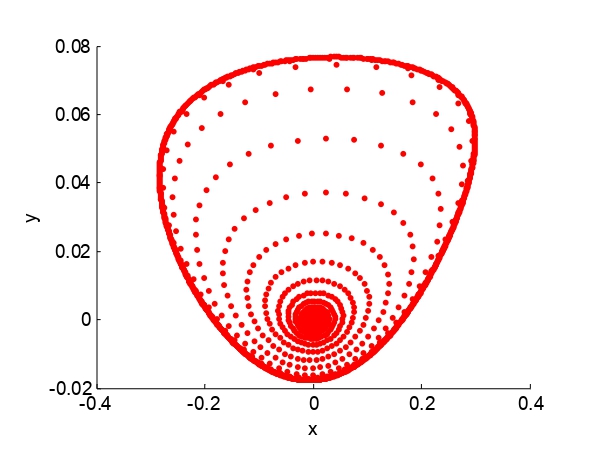}
\end{center}
\caption{Dynamics of the map \eqref{ec1} for $a$ in a neighborhood of $%
a_{NS}.$ The parameters are $m=0.02,$ $s=1.1$ and a) $a=2$ (left), b) $%
a=a_{NS}=2.0852$ and c) $a=2.1$ (right). The closed stable curve exists for $%
a>a_{NS}.$ }
\label{f3}
\end{figure}

In Fig.\ref{f3} one can see the dynamics of the map \eqref{ec1} for $a$ in a
neighborhood of $a_{NS},$ pointing out the NS bifurcation leading to a
closed stable invariant curve around the origin. Notice that, if $%
m_{0}+4\rightarrow 0,$ $m$ can be made sufficiently small.

The orbits on the invariant curve $\Gamma _{a}$ can be periodic or not \cite%
{kuz}. More exactly, if $\frac{\tilde{\varphi}-\varphi }{2\pi }=\frac{p}{q},$
with $p$ and $q\neq 0$ integers, then any orbit of (\ref{pol}) starting on
the closed invariant curve $\Gamma _{a}$ is periodic, while if $\frac{\tilde{%
\varphi}-\varphi }{2\pi }$ is irrational, no orbit from $\Gamma _{a}$ is
periodic and all orbits on $\Gamma _{a}$ are dense in $\Gamma _{a}.$

\begin{remark}
The numerical simulations that we present in this section and in the next one are
for $m$ small, since this is the case of practical relevance.
\end{remark}

\begin{figure}[h]
\begin{center}
\includegraphics[width=80mm, height=50mm]{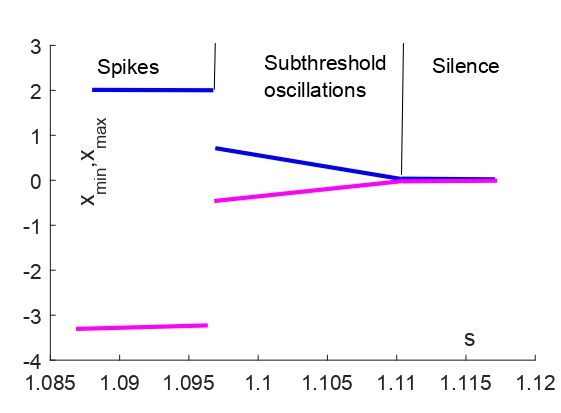}
\end{center}
\caption{Transition from silence to spiking activity through subthreshold
oscillations as $s$ decreases in the system (\protect\ref{ec0a}) for $m$
small, $m=0.02,$ and $a=2.1$}
\label{f2}
\end{figure}

\begin{figure}[h]
\begin{center}
\includegraphics[width=55mm, height=50mm]{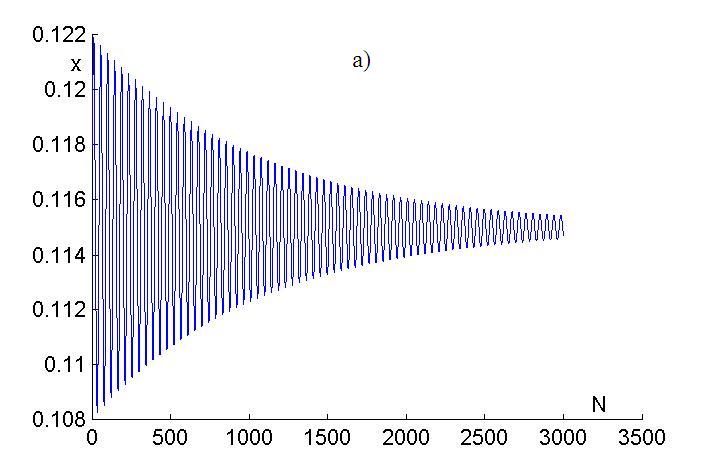} %
\includegraphics[width=55mm, height=50mm]{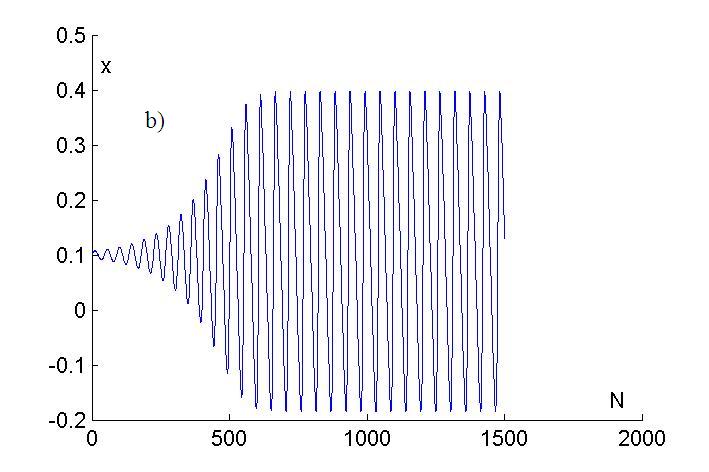} %
\includegraphics[width=55mm, height=50mm]{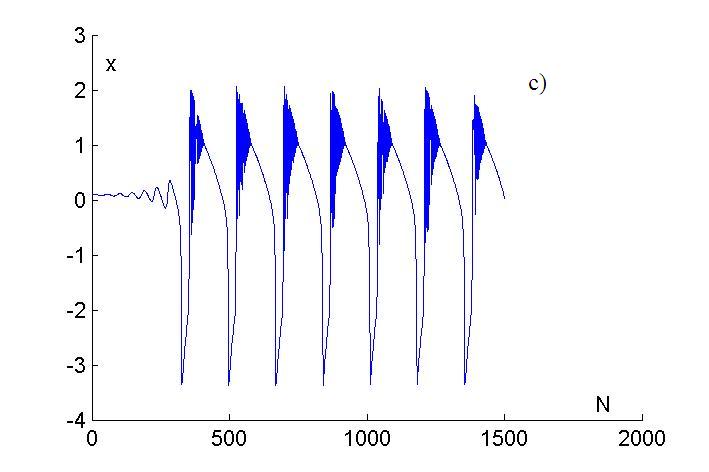}
\end{center}
\caption{Transition of the dynamics of the map \eqref{ec0a} from silence to
tonic spiking through subthreshold oscillations for $m$ small. The
parameters are $a=2.1,$ $m=0.02$ and a) $s=1.115$, b) $s=1.1$ and c) $%
s=1.09\ .$ }
\label{f2a}
\end{figure}

We illustrate in Fig.\ref{f2} and Fig.\ref{f2a} the transition from silence to spiking regime
through subthreshold oscillations. The two branches (top and bottom)
in Fig.\ref{f2} correspond to the highest and lowest values of the $x-$variable in the
initial system \eqref{ec0a}, with $f$ given by \eqref{ec0b}, for a given
value of the parameter $s;$ the other parameters are fixed at $m=0.02$ and $%
a=2.1.$ We notice that, subthreshold oscillations in Fig.\ref{f2} occur when 
$s$ ranges approximately in the interval $1.096\leq s\leq 1.11,$ which
corresponds to $2.08\leq a_{NS}\leq 2.\allowbreak 09.$ Thus, $a>a_{NS}$ in
all these cases, which implies that subthreshold periodic orbits generated
by the map \eqref{ec1} may arise in the behavior of the map (\ref{ec0a}).

The second nondegeneracy condition $e^{ik\theta _{0}}\neq 1$ leading to $%
m\in \left( 0,4\right) \smallsetminus \left\{ 2,3\right\} $ is important
because if $e^{ik\theta _{0}}=1$ for $k=1,2,3,4,$ the closed curve $\Gamma
_{a}$ may not appear at all or there might be more closed invariant curves $%
\Gamma _{a}$ bifurcating from $a=a_{NS}.$

\section{Subthreshold oscillations in the presence of noise}

\begin{figure}[h!]
\begin{center}
\includegraphics[width=55mm, height=50mm]{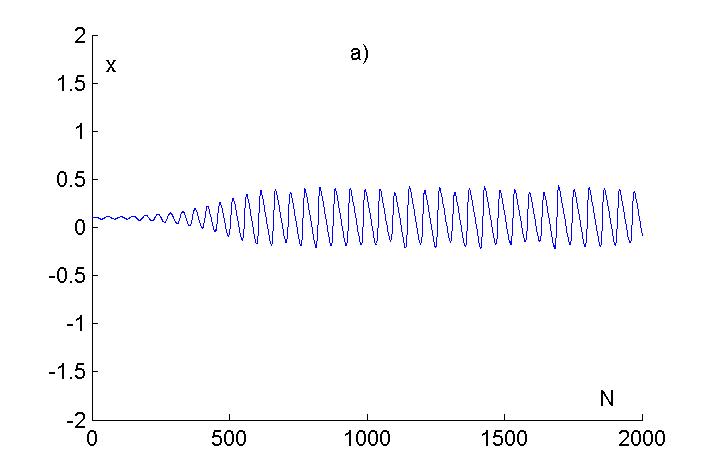} %
\includegraphics[width=55mm, height=50mm]{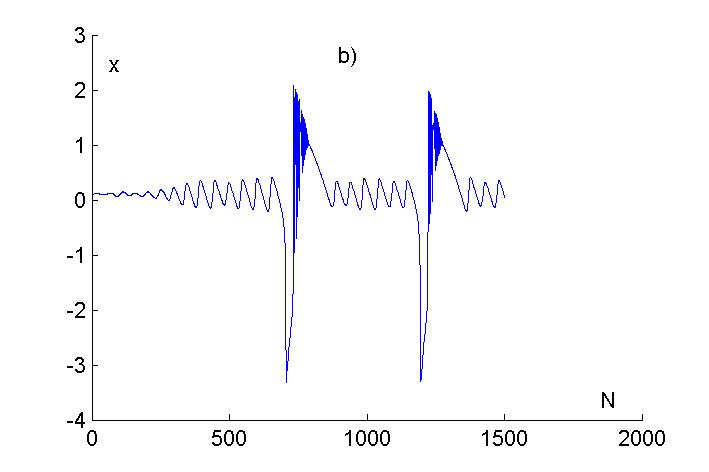} %
\includegraphics[width=55mm, height=50mm]{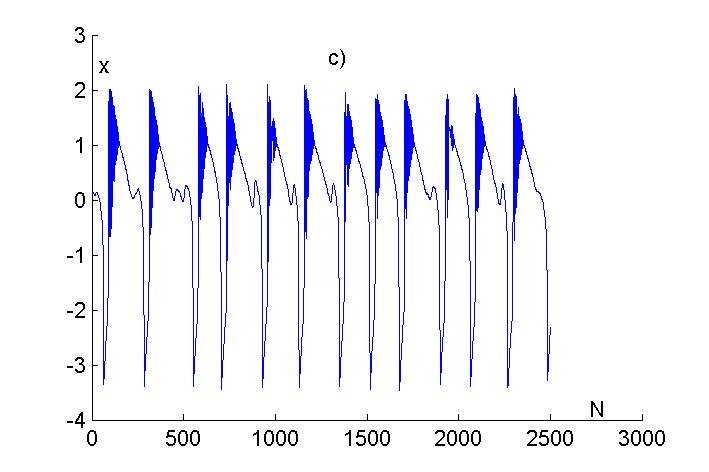}
\end{center}
\caption{Noise applied to subthreshold oscillations in the map \eqref{ec0a}
for $m$ small gives rise to spikes. The parameters are $a=2.1,$ $m=0.02,$ $%
s=1.1$ and a) $\protect\sigma=0.0001$, b) $\protect\sigma=0.0004$ and c) $%
\protect\sigma=0.004.$}
\label{f5a}
\end{figure}

Recent investigations pointed out that subthreshold activity in neurons in a
noisy medium leads to the birth of spikes \cite{bash}, \cite{bel}, \cite%
{luc19}, \cite{luc14}, \cite{luc15}, \cite{zm}. To show this property in our
model we need to apply noise to the map \eqref{ec0a} when the parameters are
in the range of subthreshold oscillations. Such ranges exist as it can be
observed numerically from the bifurcation diagrams given in Fig. \ref{f2}.
We apply noise to the slow system in the form of a Gaussian random variable
but it can be equally applied to the fast system or to both of them. The
system \eqref{ec0a} in the presence of noise added to the slow subsystem
reads: 
\begin{eqnarray}
\tilde{x} &=&f_{a}\left( x,y\right)   \label{esto1} \\
\tilde{y} &=&y-m\left( x+1-s\right) +\xi ,  \notag
\end{eqnarray}%
where $\xi $ is a Gaussian random variable of zero mean value and $\sigma $
standard deviation value. Fig. \ref{f5a} illustrates the dynamics of the
initial map \eqref{ec0a} corresponding to subthreshold oscillations in the
presence of noise. The magnitude of the noise is quantified by the standard
deviation $\sigma .$ In the case of $m$ small, if the noise is small enough,
say $\sigma =0.0001,$ no influence of it is observed in the dynamics of the
map, Fig. \ref{f5a} a). Increasing slightly $4$ times the magnitude of the
noise, $\sigma =0.0004,$ slow bursting of spikes are remarked, Fig. \ref{f5a}
b), which means that the noise is sufficient to give a start to action
potentials. A further increasing of $\sigma ,$ Fig. \ref{f5a} c), gives rise
to a tonic spiking activity.

\begin{figure}[h!]
\begin{center}
\includegraphics[width=55mm, height=50mm]{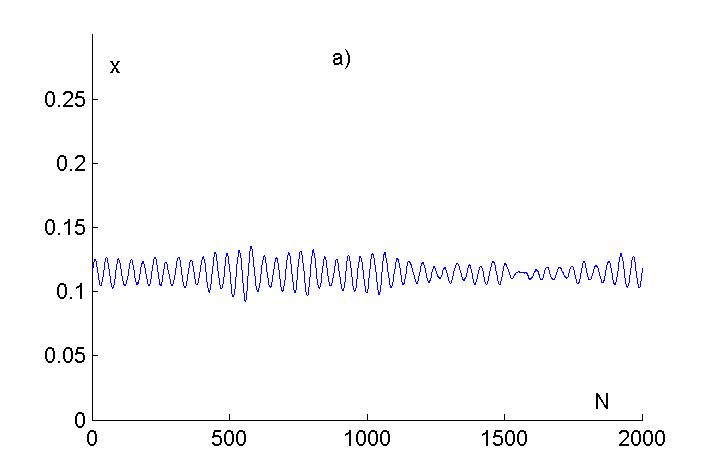} %
\includegraphics[width=55mm, height=50mm]{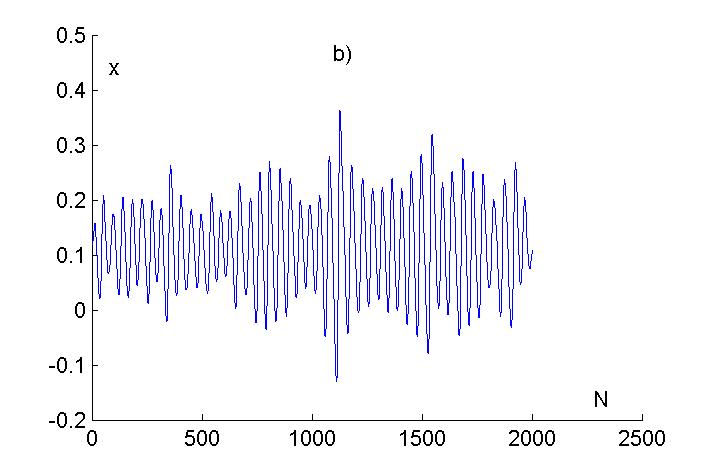} %
\includegraphics[width=55mm, height=50mm]{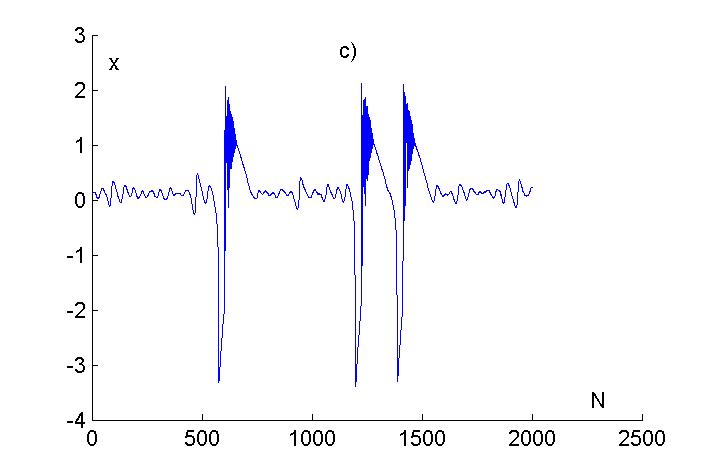}
\end{center}
\caption{Noise applied to the map \eqref{ec0a} in the silence regime for $m$
small gives rise to spikes. The parameters are $a=2.1,$ $m=0.02,$ $s=1.115$
and a) $\protect\sigma=0.0001$, b) $\protect\sigma=0.001$ and c) $\protect%
\sigma=0.002.$}
\label{f6}
\end{figure}

\begin{remark}
We notice that our map produces spikes in a noisy medium not only from the
regime of subthreshold oscillations but also from the regime of silence, as
it can be seen from Fig. \ref{f6}. Using the bifurcation diagram sketched in
Fig. \ref{f2}, we set the parameters at $a=2.1,$ $m=0.02$ and $s=1.115$
corresponding to the silence regime. While for $\sigma <0.001$ small, Figs. %
\ref{f6} a)-b), no spike is observed, a slight increasing in $\sigma ,$ $%
\sigma \geq 0.002,$ bursts of spikes are obtained in the map's dynamics,
Fig. \ref{f6} c).
\end{remark}

\section{Conclusions}

In the present work we have proposed and investigated a neuronal non-smooth
map-based model to replicate individual dynamics of a neuron. This model
comes as a theoretical tool in understanding subthreshold oscillations in
neurons. We searched for a model with a diverse dynamics in the zone of
subthreshold oscillations to describe the potential behavior of a neuron. We
studied the map particularly for $m$ small ($0<m\ll 1$) and obtained
subthreshold oscillations leading and shaping spiking activity in the map.
We sketched the transition from silence regime to spiking activity through
subthreshold oscillations in a bifurcation diagram and illustrated its
validity by numerical simulations containing waveforms of the map.

We studied then numerically if our model captures the property of action
potentials of being vulnerable to noises and found that the model has this
property. We applied noise to our map in the form of a Gaussian random
variable and obtained bursting-spiking activity. Moreover, this property of
bursts of spikes in a noisy medium is observed even though the map is in a
silence regime. This is interesting and leads us to the following
observation: if the noise is powerful enough a neuron is brought out from
its silence suddenly, without observable delays in the subthreshold phase.
This property may appear in our model due to the exponential branch of the function $f_a.$

\section{Acknowledgments}

This research was partially supported by Horizon2020-2017-RISE-777911
project.

\end{document}